\documentclass[12pt,reqno]{amsart}

\usepackage{mathrsfs,amsfonts,amsthm,amsmath,amssymb,thmtools}

\usepackage{enumitem}                
\usepackage{lineno}                 

\usepackage[hyperindex, pdfencoding = auto, psdextra, bookmarksdepth = 4]{hyperref}  
\hypersetup{
    bookmarksopen      = true,
    bookmarksopenlevel = 1,
    bookmarksnumbered  = true,
    pdftitle    = {},
    pdfauthor   = {},
    linktoc     = page,
    anchorcolor = green,
    colorlinks  = true,        
    linkcolor   = black,
    citecolor   = blue,
    filecolor   = blue,
    urlcolor    = magenta,
}

\usepackage[nameinlink]{cleveref}    

\theoremstyle{plain}

\newtheorem{theorem}{Theorem}[section]
\newtheorem{lemma}[theorem]{Lemma}
\newtheorem{proposition}[theorem]{Proposition}
\newtheorem{corollary}[theorem]{Corollary}
\newtheorem{question}[theorem]{Question}

\theoremstyle{definition}   

\newtheorem{remark}[theorem]{Remark}

\makeatletter
    \def\subsection{\@startsection{subsection}{2}%
    \z@{.5\linespacing\@plus.7\linespacing}{.3\linespacing}%
    {\normalfont\bfseries}}
\makeatother 

\makeatletter                                         
    \newcommand{\LeftEqNo}{\let\veqno\@@leqno}        
\makeatother

\numberwithin{equation}{section}                      

\begin{document}

\

\vspace{-2cm}

\title[Products of irreducible operators in separable factors]{Products of irreducible operators in factors }

\author{Minghui Ma}
\address{Minghui Ma, School of Mathematical Sciences, Dalian University of Technology, Dalian, 116024, China}
\email{minghuima@dlut.edu.cn}

\author{Junhao Shen}
\address{Junhao Shen, Department of Mathematics \& Statistics, University of New Hampshire, Durham, 03824, US}
\email{junhao.shen@unh.edu}

\author{Rui Shi}
\address{Rui Shi, School of Mathematical Sciences, Dalian University of Technology, Dalian, 116024, China}
\email{ruishi@dlut.edu.cn}
\thanks{Rui Shi, Minghui Ma, and Tianze Wang were supported by NSFC (No.12271074). 
Minghui Ma was also supported by the China Postdoctoral Science Foundation (No.2025M783077) and the Postdoctoral Fellowship Program of CPSF (No.GZC20252022).}

\author{Tianze Wang}
\address{Tianze Wang, School of Mathematical Sciences, Dalian University of Technology, Dalian, 116024, China}
\email{swan0108@mail.dlut.edu.cn}



\begin{abstract}
    Let $\mathcal{M}$ be a separable factor.
    An operator $T$ in $\mathcal{M}$ is said to be irreducible in $\mathcal{M}$ if the von Neumann algebra $W^*(T)$ generated by $T$ is an irreducible subfactor of $\mathcal{M}$, i.e., $W^*(T)'\cap\mathcal{M}=\mathbb{C}I$.
    In this paper, we show that every operator in a separable factor $\mathcal{M}$ is the product of two irreducible operators in $\mathcal{M}$, except the zero operator in factors of type $\mathrm{I}_{2n+1}$ for $n\geqslant 1$.
    This may be viewed as a multiplicative analogue of Radjavi's result which asserts that every operator on a separable Hilbert space is the sum of two irreducible operators.
\end{abstract}

\maketitle

\section{Introduction}

In 1968, P. Halmos \cite{Hal} initiated the study of irreducible operators 
in a topological sense by proving that the set of irreducible operators is a norm-dense $G_\delta$ subset of $\mathcal{B}(\mathcal{H})$, where $\mathcal{H}$ is a separable complex Hilbert space. An operator $T$ in $\mathcal{B}(\mathcal{H})$ is \emph{irreducible} if there is no nontrivial projection $P$ commuting with $T$. Note that such operators can be regarded as the fundamental building blocks to construct operators in $\mathcal{B}(\mathcal{H})$, and they help to illuminate the intrinsic structure of  $\mathcal{B}(\mathcal{H})$.
A short proof of Halmos' result \cite[Main theorem]{Hal} was provided by H. Radjavi and P. Rosenthal \cite{RR}.
Moreover, the set of irreducible operators is dense with respect to the Schatten $p$-norm for every $p>1$ (see \cite[Lemma 4.33]{Her89-book}).

Compared with the topological results of irreducible operators mentioned above, P. Fillmore and D. Topping \cite{FT} showed that irreducible operators are abundant in an algebraic sense, namely, every operator in $\mathcal{B}(\mathcal{H})$ is the sum of four irreducible operators.
H. Radjavi \cite{Rad} improved this result by proving that every operator
in $\mathcal{B}(\mathcal{H})$ is the sum of two irreducible operators. For more results about irreducible operators, the reader is referred to \cite{JW98-book,JW06}. 
Since multiplication is also an elementary algebraic operation in $\mathcal{B}(\mathcal{H})$, a natural question arises: is every operator in $\mathcal{B}(\mathcal{H})$ the product of finitely many irreducible operators?

Recall that a {\em von Neumann algebra} $\mathcal{M}$ is a self-adjoint subalgebra of $\mathcal{B}(\mathcal{H})$ that is closed in the weak-operator topology and contains $I$, the identity operator of $\mathcal{B}(\mathcal{H})$. The commutant $\mathcal{M}^{\prime}$ of $\mathcal{M}$ is the set $\{X \in \mathcal{B}(\mathcal{H}): X Y = Y X \text{ for all } Y \in \mathcal{M}\}$, which is also a von Neumann algebra. Moreover, if the center $\mathcal{Z}(\mathcal{M}):= \mathcal{M} \cap \mathcal{M}^{\prime}$ of $\mathcal{M}$ consists of scalar multiples of the identity, then $\mathcal{M}$ is a {\em factor}.
A von Neumann algebra is said to be {\em separable} if it has a separable predual space; equivalently, it acts on a separable Hilbert space.
By the type decomposition theorem for von Neumann algebras \cite[Theorem 6.5.2]{KR2}, factors can be classified into type $\mathrm{I}_n$, type $\mathrm{I}_\infty$, type $\mathrm{II}_1$, type $\mathrm{II}_\infty$, and type $\mathrm{III}$.
Note that $\mathcal{B}(\mathcal{H})$ is a factor of type $\mathrm{I}$.
As a generalization of Halmos' result \cite[Main theorem]{Hal}, the authors of \cite{WFS} proved the norm-density of irreducible operators in separable factors.
Moreover, the authors of \cite{SS20} proved that every operator in separable factors is the sum of two irreducible operators, which is a generalization of Radjavi's theorem. Inspired by recent research on operator-theoretic results in factors \cite{CFY21, CFY24, HMS26, SS19, SS20, WFS}, we study the products of irreducible operators in separable factors in this paper.

Throughout, let $\mathcal{M}$ be a separable factor.
As an analogue of irreducible operators in $\mathcal{B}(\mathcal{H})$, an operator $T$ in $\mathcal{M}$ is  {\em irreducible} (relative to $\mathcal{M}$) if the von Neumann algebra $W^*(T)$ generated by $T$ is an irreducible subfactor of $\mathcal{M}$, i.e., $W^*(T)'\cap\mathcal{M}=\mathbb{C}I$.
Denote by $\mathrm{IR}(\mathcal{M})$ the set of all irreducible operators in $\mathcal{M}$.
For every integer $n\geqslant 1$, let $\mathrm{IR}(\mathcal{M})^n$ be the set of all operators that can be decomposed as the product of $n$ irreducible operators in $\mathcal{M}$.
Since every irreducible operator in $\mathcal{M}$ is the product of two irreducible operators by \Cref{lem IR=IR2}, it follows that
\begin{equation*}
    \mathrm{IR}(\mathcal{M})\subseteq\mathrm{IR}(\mathcal{M})^2
    \subseteq\cdots
    \subseteq\mathrm{IR}(\mathcal{M})^n\subseteq\cdots
    \subseteq\bigcup_{n=1}^{\infty}\mathrm{IR}(\mathcal{M})^n\subseteq\mathcal{M}.
\end{equation*}
In the current paper, we provide a constructive proof that every nonzero operator in a separable factor $\mathcal{M}$ can be written as the product of two irreducible operators, i.e., $\mathrm{IR}(\mathcal{M})^2\setminus\{0\}=\mathcal{M}\setminus\{0\}$.
More precisely, the main results of this paper are summarized below:
\begin{enumerate} [label= $(\arabic*)$, ref= \arabic*, labelsep = 0.5 em]
    \item  If $\mathcal{M}$ is of type $\mathrm{I}_n$, i.e., $\mathcal{M}\cong M_n(\mathbb{C})$, then every operator in $\mathcal{M}$ is the product of two irreducible operators, except the zero operator in factors of type $\mathrm{I}_{2k+1}$ for $k\geqslant 1$ (\Cref{prop 0}, \Cref{thm I-III}).
    In other words, we have
    \begin{equation*}
        \mathrm{IR}(\mathcal{M})^2=
        \begin{cases}
            \mathcal{M}, & \mbox{if}~n=1,2,4,6,8,\ldots, \\
            \mathcal{M}\setminus\{0\}, & \mbox{otherwise}.
        \end{cases}
    \end{equation*}
    Moreover, $0\in\mathrm{IR}(\mathcal{M})^3$ (\Cref{rem 0}).
    \item  If $\mathcal{M}$ is infinite dimensional, then every operator in $\mathcal{M}$ is the product of two irreducible operators (\Cref{prop 0}, \Cref{thm I-III}, \Cref{thm II-infty}, and \Cref{thm II-1}).
    In other words, we have $\mathrm{IR}(\mathcal{M})^2=\mathcal{M}$.
\end{enumerate}

Recall that an operator in $\mathcal{B}(\mathcal{H})$ is said to be {\em strongly irreducible} if $XTX^{-1}$ is irreducible for each invertible operator $X$ in $\mathcal{B}(\mathcal{H})$, which are a natural analogue of Jordan blocks in $M_n(\mathbb{C})$.
Denote by $\mathrm{SIR}(\mathcal{H})$ the set of strongly irreducible operators in $\mathcal{B}(\mathcal{H})$.
Compared with Radjavi's result \cite{Rad}, C. Jiang and P. Wu \cite[Theorem 2.1]{JW98} proved that every operator in $\mathcal{B}(\mathcal{H})$ is the sum of three strongly irreducible operators.
It follows from \cite[Theorem 3.1.1]{JW06} that every operator in $\mathcal{B}(\mathcal{H})$ is the sum of two strongly irreducible operators.
Note that the rank of every Jordan block in $M_n(\mathbb{C})$ is at least $n-1$.
As a consequence, the rank of the product of two strongly irreducible operators in $M_n(\mathbb{C})$ is at least $n-2$.
Therefore, we cannot expect every operator to be the product of two strongly irreducible operators.
We end the introduction with the following question.

\begin{question}\label{que SIR}
    Let $\mathcal{H}$ be a separable complex Hilbert space. Is every  operator in $\mathcal{B}(\mathcal{H})$ the product of finitely many strongly irreducible operators? In other words, is it true that $\bigcup_{n=1}^{\infty}\mathrm{SIR}(\mathcal{H})^n=\mathcal{B}(\mathcal{H})$?     
\end{question}



\section{Preliminaries}  \label{section-2}

Throughout, let $\mathcal{M}$ be a separable factor.
Recall that 
$\mathrm{IR}(\mathcal{M})$ is a norm-dense $G_\delta$ subset of $\mathcal{M}$ by \cite[Theorem 2.1]{WFS}.
In particular, there are invertible irreducible operators in $\mathcal{M}$.
This fact will be used frequently in the sequel.
The following lemma shows that $\mathrm{IR}(\mathcal{M})$ is a subset of $\mathrm{IR}(\mathcal{M})^2$.
As usual, we denote by $\sigma(T)$ the spectrum of an operator $T$ in $\mathcal{M}$.

\begin{lemma}\label{lem IR=IR2}
    Every irreducible operator in a separable factor $\mathcal{M}$ is the product of an irreducible operator and an invertible irreducible operator in $\mathcal{M}$.
\end{lemma}

\begin{proof}
    Let $T$ be irreducible in $\mathcal{M}$.
    We construct two operators in $\mathcal{M}$ of the form
    \begin{equation*}
        X:= T(\lambda I-T)^{-1}\quad\text{and}\quad Y:= \lambda I-T,
    \end{equation*}
    where $\lambda$ is a scalar in $\mathbb{C}\setminus\sigma(T)$.
    Then $Y$ is invertible in $\mathcal{M}$ and $T=XY$.
    A direct computation shows that $T=\lambda X(I+X)^{-1}=\lambda I-Y$.
    Thus, both $X$ and $Y$ are irreducible in $\mathcal{M}$.
\end{proof}

In the following lemma, we characterize the relative commutant for a class of operators in a separable factor.

\begin{lemma}\label{lem-IR-opt-direct-sum-commutant}
    Let $\mathcal{M}$ be a separable factor, $\{P_1,P_2\}$ nonzero orthogonal projections in $\mathcal{M}$, and $P=P_1+P_2$.
    Suppose that $T_j$ is an irreducible operator in $P_j\mathcal{M}P_j$ for $j=1,2$.
    If there exists no partial isometry $V$ in $P_1\mathcal{M}P_2$ satisfying $V^*V=P_2$, $VV^*=P_1$, and $ T_1 V = V T_2$, then $W^*(T_1+T_2)'\cap P\mathcal{M}P=\mathbb{C}P_1\oplus\mathbb{C}P_2$.
\end{lemma}

\begin{proof}
    Let $Q$ be a projection in $P\mathcal{M}P$ commuting with $T_1+T_2$.
    For simplicity, let $Q_{ij}=P_iQP_j$ for $1\leqslant i,j\leqslant 2$.
    Then $T_1Q_{12}=Q_{12}T_2$ and $T_2Q_{21}=Q_{21}T_1$.
    It follows that
    \begin{equation*}
        T_1Q_{12}Q_{21}=Q_{12}Q_{21}T_1\quad\text{and}\quad
        T_2Q_{21}Q_{12}=Q_{21}Q_{12}T_2.
    \end{equation*}
    By the irreducibility of $T_1$ and $T_2$, we have
    \begin{equation*}
        Q_{12}Q_{21}=c^2P_1\quad\text{and}\quad Q_{21}Q_{12}=c^2P_2,
    \end{equation*}
    where $c=\|Q_{12}\|$.
    We claim that $Q_{12}=0$.
    Otherwise, let
    \begin{equation*}
        E_{11}=P_1,\quad E_{12}=c^{-1}Q_{12},\quad E_{21}=c^{-1}Q_{21},
        \quad\text{and}\quad E_{22}=P_2.
    \end{equation*}
    Then $\{E_{ij}\}_{i,j=1}^2$ is a system of matrix units in $P\mathcal{M}P$ such that $T_1E_{12}=E_{12}T_2$.
    That contradicts the assumption. Thus, $Q = Q_{11} + Q_{22}$. By the irreducibility of $T_1$ and $T_2$, we can complete the proof.
\end{proof}

\begin{remark}\label{lem IR-sum-commutant}
    For $T_1$ and $T_2$ in \Cref{lem-IR-opt-direct-sum-commutant}, notice that if $\|T_1\|\ne\|T_2\|$, then there is no partial isometry $V$ in $P_1\mathcal{M}P_2$ satisfying the conditions in \Cref{lem-IR-opt-direct-sum-commutant}. 
    We will apply a combination of this fact and  \Cref{lem-IR-opt-direct-sum-commutant} in the proof of \Cref{prop key}.
\end{remark}

Let $\mathcal{A}$ be a unital Banach algebra.
Recall that for any $A,B\in\mathcal{A}$, the {\em Rosenblum operator $\mathcal{T}_{AB}$} on $\mathcal{A}$ is defined by $\mathcal{T}_{AB}(X):= AX-XB$.
By \cite[Corollary 3.1]{Ros}, the operator $\mathcal{T}_{AB}$ is invertible whenever $\sigma(A)\cap\sigma(B)= \varnothing$; see also \cite[Corollary 0.13]{RR73}.
We prove an analogue for Banach bimodules in \Cref{lem rosenblum}.
The case $\mathcal{A}=P_1\mathcal{M}P_1$, $\mathcal{B}=P_2\mathcal{M}P_2$, and $\mathcal{X}=P_1\mathcal{M}P_2$ will be directly used in the proof of \Cref{lem IR-triangular}.

\begin{lemma}\label{lem rosenblum}
    Let $\mathcal{X}$ be a Banach $\mathcal{A}$-$\mathcal{B}$-bimodule, where $\mathcal{A}$ and $\mathcal{B}$ are unital Banach algebras.
    For any $A\in\mathcal{A}$ and $B\in\mathcal{B}$, if $\sigma_{\mathcal{A}}(A) \cap \sigma_{\mathcal{B}}(B) = \varnothing$, then the Rosenblum operator $\mathcal{T}_{AB}$ on $\mathcal{X}$  is invertible, where $\mathcal{T}_{AB}$ is defined by $\mathcal{T}_{AB}(X) := AX-XB$ for every $X$ in $\mathcal{X}$.
\end{lemma}

\begin{proof}
    By assumption, there exists a Cauchy domain $\Omega$ such that $\sigma_{\mathcal{B}}(B)\subseteq\Omega$ and $\sigma_{\mathcal{A}}(A)\cap\overline{\Omega}= \varnothing$.
    For every $z\in\partial\Omega$, the elements $A-zI_{\mathcal{A}}$ and $zI_{\mathcal{B}}-B$ are invertible in $\mathcal{A}$ and $\mathcal{B}$, respectively.
    Since $(A-zI_{\mathcal{A}})X+X(zI_{\mathcal{B}}-B)=\mathcal{T}_{AB}(X)$, we have
    \begin{equation*}
        X(zI_{\mathcal{B}}-B)^{-1}+(A-zI_{\mathcal{A}})^{-1}X
        =(A-zI_{\mathcal{A}})^{-1}\mathcal{T}_{AB}(X)(zI_{\mathcal{B}}-B)^{-1}.
    \end{equation*}
    It follows that
    \begin{equation*}
    X=\int_{\partial\Omega}
    (A-zI_{\mathcal{A}})^{-1}\mathcal{T}_{AB}(X)(zI_{\mathcal{B}}-B)^{-1}dz.
    \end{equation*}
    Therefore, $\mathcal{T}_{AB}$ is invertible and
    \begin{equation*}
        \mathcal{T}_{AB}^{-1}(Y)=\int_{\partial\Omega}
        (A-zI_{\mathcal{A}})^{-1}Y(zI_{\mathcal{B}}-B)^{-1}dz.
    \end{equation*}
    This completes the proof.
\end{proof}

For every operator $T$ in $\mathcal{M}$, its range projection is denoted by $R(T)$.
We provide a useful criterion for the irreducibility of lower-triangular operators as follows.

\begin{lemma}\label{lem IR-triangular}
    Let $P_1$ and $P_2$ be projections in a separable factor $\mathcal{M}$ with $I=P_1+P_2$.
    Suppose that $T_{11}$ is an irreducible operator in $P_1\mathcal{M}P_1$, $T_{21}$ is an operator in $P_2\mathcal{M}P_1$ with range projection $R(T_{21})=P_2$, and $T_{22}$ is an operator in $P_2\mathcal{M}P_2$ such that
    \begin{equation*}
        \sigma_{P_1\mathcal{M}P_1}(T_{11})\cap
        \sigma_{P_2\mathcal{M}P_2}(T_{22})= \varnothing.
    \end{equation*}
    Let $T=T_{11}+T_{21}+T_{22}$.
    Then $T$ is irreducible in $\mathcal{M}$.
\end{lemma}

\begin{proof}
    Let $Q$ be a projection in $\mathcal{M}$ commuting with $T$.
    We prove that $Q$ is trivial.
    For simplicity, write $Q_{ij} := P_iQP_j$ for $1\leqslant i,j\leqslant 2$.
    Then $T_{11}Q_{12}=Q_{12}T_{22}$ and $Q_{12}$ lies in $P_1\mathcal{M}P_2$.
    It follows from \Cref{lem rosenblum} that $Q_{12}=0$, and hence $Q_{21}=Q_{12}^*=0$.
    Consequently, $Q_{11}$ is a projection in $P_1\mathcal{M}P_1$ commuting with $T_{11}$.
    By the irreducibility of $T_{11}$, we have $Q_{11}=0$ or $P_1$.
    Without loss of generality, we may assume that $Q_{11}=0$.
    Then $Q_{22}T_{21}=0$.
    Since $R(T_{21})=P_2$, we obtain that $Q_{22}=0$, i.e., $Q=0$.
    Therefore, $T$ is irreducible in $\mathcal{M}$.
\end{proof}

Recall that two projections $P$ and $Q$ in $\mathcal{M}$ are said to be (\emph{Murray-von Neumann}) {\em  equivalent}, denoted by $P\sim Q$, if there exists a partial isometry $V$ in $\mathcal{M}$ such that $V^*V=P$ and $VV^*=Q$.
If $Q$ is equivalent to a subprojection of $P$ in $\mathcal{M}$, then we say that $Q$ is weaker than $P$, denoted by $Q\precsim P$ or $P\succsim Q$.
As an application of \Cref{lem IR-triangular}, we show that a class of operators can be written as the product of two irreducible operators.

\begin{lemma}\label{lem half-exceed}
    Let $\mathcal{M}$ be a separable factor.
    Suppose that $T$ is an operator in $\mathcal{M}$ and $P$ is a projection in $\mathcal{M}$ such that
    \begin{equation*}
        PT=TP\quad\text{and}\quad P\succsim I-P.
    \end{equation*}
    If $TP$ is the product of an irreducible operator and an invertible irreducible operator in $P\mathcal{M}P$, then $T$ is the product of two irreducible operators in $\mathcal{M}$.
\end{lemma}

\begin{proof}
    Let $P_1=P$, $P_2=I-P$, and $T_j=TP_j\in P_j\mathcal{M}P_j$ for $j=1,2$.
    By assumption, we can write $T_1=AB$, where $A$ is an irreducible operator and $B$ is an invertible irreducible operator in $P_1\mathcal{M}P_1$.
    Then there exists a large scalar $\lambda$ in $\mathbb{C}\setminus\sigma_{P_1\mathcal{M}P_1}(A)$ such that $\sigma_{P_1\mathcal{M}P_1}(B)
    \cap\sigma_{P_2\mathcal{M}P_2}(\lambda^{-1}T_2)= \varnothing$.
    Since $P_1\succsim P_2$, there exists a partial isometry $V$ in $P_2\mathcal{M}P_1$ such that $V^*V\leqslant P_1$ and $VV^*=P_2$.
    Let
    \begin{equation*}
        X := A+\lambda V+\lambda P_2\quad\text{and}\quad Y := B-VB+\lambda^{-1}T_2.
    \end{equation*}
    It is clear that $T=XY$.
    By \Cref{lem IR-triangular}, both $X$ and $Y$ are irreducible in $\mathcal{M}$.
    This completes the proof.
\end{proof}

\begin{remark}\label{rem half-exceed}
    The condition on $TP$ in \Cref{lem half-exceed} is reasonable by \Cref{lem IR=IR2}.
\end{remark}

The following lemma states that if the projections $P$ and $I-P$ in \Cref{lem half-exceed} are equivalent, then we can remove the condition on $TP$.

\begin{lemma}\label{lem half-exact}
    Let $\mathcal{M}$ be a separable factor.
    Suppose that $T$ is an operator in $\mathcal{M}$ and $P$ is a projection in $\mathcal{M}$ such that
    \begin{equation*}
        PT=TP\quad\text{and}\quad P\sim I-P.
    \end{equation*}
    Then $T$ is the product of two irreducible operators in $\mathcal{M}$.
\end{lemma}

\begin{proof}
    By the assumption $P\sim I-P$, there exists a system of matrix units $\{E_{ij}\}_{i,j=1}^{2}$ in $\mathcal{M}$ such that $E_{11}=P$ and $E_{22}=I-P$.
    Since $TP=PT$, it follows that $T=T_1+T_2$, where $T_j=TE_{jj}\in E_{jj}\mathcal{M}E_{jj}$ for $j=1,2$.
    By \cite[Theorem 2.1]{WFS}, let $X_1$ be an invertible irreducible operator in $E_{11}\mathcal{M}E_{11}$ and $X_2=E_{21}X_1E_{12}\in E_{22}\mathcal{M}E_{22}$.
    The inverse of $X_j$ in $E_{jj}\mathcal{M}E_{jj}$ is denoted by $X_j^{-1}$ for $j=1,2$.
    Then there exists a nonzero small scalar $\lambda$ such that
    \begin{equation*}
        \sigma_{E_{11}\mathcal{M}E_{11}}(\lambda T_1X_1)
        \cap\sigma_{E_{22}\mathcal{M}E_{22}}(X_2)= \varnothing =
        \sigma_{E_{11}\mathcal{M}E_{11}}(\lambda^{-1}X_1^{-1})
        \cap\sigma_{E_{22}\mathcal{M}E_{22}}(X_2^{-1}T_2).
    \end{equation*}
    We construct two operators in $\mathcal{M}$ as follows
    \begin{equation*}
        X := \lambda T_1X_1+\lambda E_{21}X_1+X_2\quad\text{and}\quad
        Y := \lambda^{-1}X_1^{-1}-X_2^{-1}E_{21}+X_2^{-1}T_2.
    \end{equation*}
    Then $T=XY$, and the irreducibility of $X$ and $Y$ follows from \Cref{lem IR-triangular}.
\end{proof}

The following corollary is a direct consequence of \Cref{lem half-exact}.

\begin{corollary}  \label{cor finite-rank}
    Let $\mathcal{H}$ be a separable infinite-dimensional complex Hilbert space.
    Then every finite-rank operator in $\mathcal{B}(\mathcal{H})$ is the product of two irreducible operators in $\mathcal{B}(\mathcal{H})$.
\end{corollary}

In the following proposition, we completely determine that the zero operator is the product of two irreducible operators in certain factors.
For an operator $X$, we denote by $\operatorname{ran} X$ and $\ker X$ the range space and the null space of $X$, respectively.

\begin{proposition}\label{prop 0}
    Let $\mathcal{M}$ be a separable factor.
    Then the zero operator is the product of two irreducible operators in $\mathcal{M}$ if and only if $\mathcal{M}$ is not of type $\mathrm{I}_{2n+1}$ for each $n\geqslant 1$.
\end{proposition}

\begin{proof}
    If $\mathcal{M}$ is not of type $I_{2n-1}$ for each integer $n\geqslant 1$, then there exists a projection $P$ in $\mathcal{M}$ such that $P\sim I-P$.
    By \Cref{lem half-exact}, the zero operator is the product of two irreducible operators in $\mathcal{M}$.
    Note that every scalar is irreducible in $\mathbb{C}$ by definition.
    We claim that the zero operator is not the product of two irreducible operators in $M_{2n+1}(\mathbb{C})$.
    Otherwise, we assume that $0=XY$ for some irreducible operators $X$ and $Y$ in $M_{2n+1}(\mathbb{C})$.
    Since the nonzero projection $R(X)\vee R(X^*)$ commutes with $X$, we must have $R(X)\vee R(X^*)=I$.
    It follows that $\dim\operatorname{ran} X \geqslant n+1$, i.e., $\dim\ker X \leqslant n$.
    Similarly, $\dim\operatorname{ran} Y\geqslant n+1$.
    Thus, $\operatorname{ran} Y$ is not a subspace of $\ker X$.
    That is a contradiction.
    We complete the proof.
\end{proof}

\begin{remark}\label{rem 0}
    The zero operator in $M_{2n+1}(\mathbb{C})$ is the product of three irreducible operators.
    Actually, let $J$ be a Jordan block with eigenvalue $0$.
    Then $J$ is irreducible and there exists a rank-one projection $P$ such that $JP=0$.
    By \Cref{thm I-III}, $P$ is the product of two irreducible operators.
\end{remark}

\section{Products of irreducible operators in factors of types \texorpdfstring{$\mathrm{I}$}{I} and \texorpdfstring{$\mathrm{III}$}{III}} \label{section-3}

In this section, we prove that every nonzero operator is the product of two irreducible operators in separable factors of type $\mathrm{I}$ and type $\mathrm{III}$.
The zero operator has been handled in \Cref{prop 0}.

Recall that an operator $T$ in $\mathcal{B}(\mathcal{H})$ has {\em closed range} if its range space 
is a closed subspace of $\mathcal{H}$; equivalently, there exists a scalar $c>0$ such that $T^*T\geqslant c R(T^*)$.
In the following lemma, we show that the direct sum of certain operators with closed range can be written as the product of two irreducible operators.
This lemma will be used to derive \Cref{prop closed-range}.

\begin{lemma}\label{lem closed-range}
    Let $\mathcal{M}$ be a separable factor and $\{P_1,P_2\}$ nonzero orthogonal projections in $\mathcal{M}$.
    Suppose that $T_j$ is a nonzero operator in $P_j\mathcal{M}P_j$ and is the product of a nonzero irreducible operator and an invertible irreducible operator in $P_j\mathcal{M}P_j$ for each $j=1,2$.
    If $T_1$ has closed range, then $T_1+T_2$ is the product of a nonzero irreducible operator and an invertible irreducible operator in $P\mathcal{M}P$, where $P=P_1+P_2$.
\end{lemma}

\begin{proof}
    Let $T_j=A_jB_j$, where $A_j$ is a nonzero irreducible operator and $B_j$ is an invertible irreducible operator in $P_j\mathcal{M}P_j$ for $j=1,2$.
    Since $T_1$ has closed range and $B_1$ is invertible, $A_1$ also has closed range.
    Since $T_1=(\lambda A_1)(\lambda^{-1}B_1)$ for each nonzero scalar $\lambda$, we may assume that
    \begin{equation*}
        A_1^*A_1\geqslant(1+2\|A_2\|^2)R(A_1^*)\quad\text{and}\quad
        \sigma_{P_1\mathcal{M}P_1}(B_1)\cap\sigma_{P_2\mathcal{M}P_2}(B_2)= \varnothing.
    \end{equation*}
    Clearly, there is a partial isometry $V$ in $P_2\mathcal{M}P_1$ such that $A_2V\ne 0$.
    If $R(A_1^*)\ne P_1$, we further assume that $A_2VA_1^*=0$.
    Let
    \begin{equation*}
        X=A_1+A_2V+A_2\quad\text{and}\quad Y=B_1-VB_1+B_2.
    \end{equation*}
    Then $T=XY$.
    It is clear that $Y$ is irreducible in $P\mathcal{M}P$ by \Cref{lem IR-triangular}.
    Next, we prove the irreducibility of $X$ in $P\mathcal{M}P$ in two cases. Let $Q$ be a projection in $P\mathcal{M}P$ commuting with $X$. For simplicity, let $Q_{ij}=P_iQP_j$ for $1\leqslant i,j\leqslant 2$.

    \noindent \textbf{Case I.} $R(A_1^*)=P_1$.
    That $XQ = QX$ yields $A_1Q_{12}=Q_{12}A_2$. Hence
    \begin{equation*}
        (1+2\|A_2\|^2)Q_{12}^*Q_{12}\leqslant Q_{12}^*A_1^*A_1Q_{12}
        =A_2^*Q_{12}^*Q_{12}A_2\leqslant\|Q_{12}\|^2A_2^*A_2.
    \end{equation*}
    Taking the operator norm in both sides, we have $(1+2\|A_2\|^2)\|Q_{12}\|^2\leqslant\|A_2\|^2\|Q_{12}\|^2$.
    It follows that $Q_{12}=0$. This, combined with the irreducibility of $A_j$, implies that $Q_{jj}= 0$ or $P_j$ for $j=1,2$. Since $A_2V \neq 0$ and $Q_{22} A_2 V = A_2 V Q_{11}$, we can obtain that $Q = 0$ or $P$.
    Therefore, $X$ is irreducible in $P\mathcal{M}P$.
    
    \noindent \textbf{Case II.} $R(A_1^*)\ne P_1$.
    By the assumption $A_2VA_1^*=0$, it is routine to verify that
    \begin{equation*}
        XX^*=A_1A_1^*+A_2(P_2+VV^*)A_2^*.    
    \end{equation*}
    By the well-known equality $\sigma(AB)\setminus\{0\}=\sigma(BA)\setminus\{0\}$ for all operators $A$ and $B$, we have $\sigma_{P_1\mathcal{M}P_1}(A_1A_1^*)
    \subseteq\{0\}\cup\sigma_{P_1\mathcal{M}P_1}(A_1^*A_1)
    \subseteq\{0\}\cup[1+2\|A_2\|^2,\infty)$, and then
    \begin{equation*}
      \sigma_{P_1\mathcal{M}P_1}(A_1A_1^*)
      \cap\sigma_{P_2\mathcal{M}P_2}(A_2(P_2+VV^*)A_2^*)\subseteq\{0\}.
    \end{equation*}
    It follows that $R(A_1)\in W^*(XX^*)\subseteq W^*(X)$.
    Thus, $A_1=R(A_1)X\in W^*(X)$.
    Note that $R(A_1)\vee R(A_1^*)$ is a nonzero projection in $P_1\mathcal{M}P_1$ commuting with $A_1$.
    By the irreducibility of $A_1$, we have $P_1=R(A_1)\vee R(A_1^*)\in W^*(X)$. Thus, we have that $QP_1 = P_1Q$. This yields that $Q = Q_{11}+ Q_{22}$. Similar to the proof of \textbf{Case I}, we achieve that $Q=0$ or $P$.
    Therefore, $X$ is irreducible in $P\mathcal{M}P$.
\end{proof}

The following proposition is an immediate consequence of \Cref{lem closed-range}.

\begin{proposition}\label{prop closed-range}
    Let $\mathcal{M}$ be a separable factor and $\{P_j\}_{j=1}^{k}$ a finite orthogonal family of nonzero projections in $\mathcal{M}$.
    Suppose that $T_j$ in $P_j\mathcal{M}P_j$ is the product of a nonzero irreducible operator and an invertible irreducible operator in $P_j\mathcal{M}P_j$ for all $1\leqslant j\leqslant k$.
    
    If $T_j$ has closed range for each $1\leqslant j\leqslant k-1$, then $\sum_{j=1}^{k}T_j$ is the product of a nonzero irreducible operator and an invertible irreducible operator in $P\mathcal{M}P$, where $P=\sum_{j=1}^{k}P_j$.
\end{proposition}

\begin{proof}
    The case $k=2$ follows from \Cref{lem closed-range}.
    If $k\geqslant 3$, then by induction, $\sum_{j=1}^{k-1}T_j$ is the product of a nonzero irreducible operator and an invertible irreducible operator in $Q\mathcal{M}Q$, where $Q=\sum_{j=1}^{k-1}P_j$.
    Since $T_j$ has closed range for each $1\leqslant j\leqslant k-1$, the operator $\sum_{j=1}^{k-1}T_j$ has closed range.
    Thus, the proof is completed by \Cref{lem closed-range}.
\end{proof}

The next lemma shows that $\mathrm{IR}(P\mathcal{M}P)^2\setminus\{0\}\subseteq\mathrm{IR}(\mathcal{M})^2$ for certain projections $P$ in $\mathcal{M}$.
Note that if $\mathcal{M}=M_{2n+1}(\mathbb{C})$ for some $n\geqslant 1$ and $P$ is a minimal projection in $\mathcal{M}$, then $0\in\mathrm{IR}(P\mathcal{M}P)^2$, but $0\notin\mathrm{IR}(\mathcal{M})^2$ by \Cref{prop 0}.

\begin{lemma}  \label{lem add-0-even}
    Let $\mathcal{M}$ be a separable factor.
    Suppose that $P$, $P_1$, and $P_2$ are projections in $\mathcal{M}$ such that $I=P+P_1+P_2$ and $P_1\sim P_2$.
    Let $T$ be a nonzero operator in $P \mathcal{M} P$.
    If $T$ is the product of two irreducible operators in $P\mathcal{M}P$, then $T$ is the product of two irreducible operators in $\mathcal{M}$.
\end{lemma}

\begin{proof}
    By the assumption $P_1\sim P_2$, there is a system of matrix units $\{E_{ij}\}_{i,j=1}^{2}$ such that $E_{jj}=P_j$ for $j=1,2$.
    Let $T=AB$, where $A$ and $B$ are nonzero irreducible operators in $P\mathcal{M}P$.
    Then there is a partial isometry $V$ in $P_2\mathcal{M}P$ such that $VB\ne 0$.
    Let $X_1$ be an invertible irreducible operator in $P_1\mathcal{M}P_1$ such that
    \begin{equation*}
        \sigma_{P\mathcal{M}P}(A)\cap\sigma_{P_1\mathcal{M}P_1}(X_1)= \varnothing\quad
        \text{and}\quad
        \sigma_{P\mathcal{M}P}(B)\cap\sigma_{P_1\mathcal{M}P_1}(X_1)= \varnothing.
    \end{equation*}
    Let $X_2=E_{21}X_1E_{12}\in P_2\mathcal{M}P_2$.
    We define
    \begin{equation*}
        X:= A+X_2V+E_{21}+X_2\quad\text{and}\quad
        Y:= B+X_1-VB-E_{21}.
    \end{equation*}
    Write $P_0=P$.
    For the reader's convenience, we write $X$ and $Y$ in matrix form 
    \begin{equation*}
        X=
        \begin{pmatrix}
            A & 0 & 0 \\
            0 & 0 & 0 \\
            X_2V & E_{21} & X_2
        \end{pmatrix}
        \begin{matrix}
            \operatorname{ran} P_0 \\
            \operatorname{ran} P_1 \\
            \operatorname{ran} P_2
        \end{matrix}\quad\text{and}\quad
        Y=
        \begin{pmatrix}
            B & 0 & 0 \\
            0 & X_1 & 0 \\
            -VB & -E_{21} & 0
        \end{pmatrix}
        \begin{matrix}
            \operatorname{ran} P_0 \\
            \operatorname{ran} P_1 \\
            \operatorname{ran} P_2
        \end{matrix}.    
    \end{equation*}
    A direct computation shows that $T=XY$.
    We claim that both $X$ and $Y$ are irreducible in $\mathcal{M}$ as in \Cref{lem IR-triangular}.

    Let $Q$ be a projection in $\mathcal{M}$.
    For simplicity, let $Q_{ij}=P_iQP_j$ for $0\leqslant i,j\leqslant 2$.
    If $Q$ commutes with $X$, then $AQ_{02}=Q_{02}X_2$ and $Q_{12}X_2=0$.
    It follows that $Q_{02}=0$ and $Q_{12}=0$ by \Cref{lem rosenblum} and $\sigma_{P\mathcal{M}P}(A)\cap\sigma_{P_2\mathcal{M}P_2}(X_2)= \varnothing$.
    Consequently, $Q_{22}$ is a projection in $P_2\mathcal{M}P_2$ commuting with $X_2$.
    By the irreducibility of $X_2$, we may assume that $Q_{22}=0$.
    Then $AQ_{01}=0$ and $Q_{10}A=0$, i.e., $Q_{10}A^*=Q_{10}A=0$.
    Since $R(A)\vee R(A^*)$ is a nonzero projection in $P_0\mathcal{M}P_0$ commutes with $A$, we have $P_0=R(A)\vee R(A^*)$.
    Hence $Q_{10}=0$.
    Thus, $AQ_{00}=Q_{00}A$, $X_2VQ_{00}=0$, and $E_{21}Q_{11}=0$.
    It is clear that $Q_{11}=0$.
    By the irreducibility of $A$, we see that $Q_{00}=0$.
    Therefore, $Q=0$, and hence $X$ is irreducible in $\mathcal{M}$.
    
    If $Q$ commutes with $Y$, then $X_1Q_{12}=0$ and hence $Q_{12}=0$.
    It follows that $X_1Q_{10}=Q_{10}B$.
    Hence $Q_{10}=0$ by \Cref{lem rosenblum} and $\sigma_{P\mathcal{M}P}(B)\cap\sigma_{P_1\mathcal{M}P_1}(X_1)= \varnothing$.
    Consequently, $Q_{11}$ is a projection in $P_1\mathcal{M}P_1$ commuting with $X_1$.
    By the irreducibility of $X_1$, we may assume that $Q_{11}=0$.
    Then $Q_{02}E_{21}=0$ and $Q_{22}E_{21}=0$.
    Hence $Q_{02}=0$ and $Q_{22}=0$.
    Thus, $Q_{00}$ commutes with $B$ and then $Q_{00}=0$ or $P_0$.
    Since $VBQ_{00}=0$ and $VB\ne 0$, we have $Q_{00}=0$.
    Therefore, $Y$ is also irreducible in $\mathcal{M}$.
    This completes the proof.
\end{proof}

\begin{corollary}  \label{cor add-0s}
    Let $\mathcal{M}$ be a separable factor, $P$ a nonzero projection in $\mathcal{M}$, and $T$ in $P\mathcal{M}P$ the product of a nonzero irreducible operator and an invertible irreducible operator in $P\mathcal{M}P$. Then $T$ is the product of two irreducible operators in $\mathcal{M}$.
\end{corollary}

\begin{proof}
    If there are projections $P_1$ and $P_2$ in $\mathcal{M}$ such that $I-P=P_1+P_2$ and $P_1\sim P_2$, then we complete the proof by \Cref{lem add-0-even}.
    Otherwise, $(I-P)\mathcal{M}(I-P)$ is a factor of type $\mathrm{I}_{2n-1}$ for an integer $n\geqslant 1$ and $\mathcal{M}$ is of type $\mathrm{I}$.
    Let $P_1$ be a minimal projection in $(I-P)\mathcal{M}(I-P)$.
    Then $P_1\precsim P$.
    By \Cref{lem half-exceed}, $T$ is the product of two irreducible operators in $(P+P_1)\mathcal{M}(P+P_1)$.
    Then, we complete the proof by \Cref{lem add-0-even}.
\end{proof}

We are now ready to prove the main theorem in this section.

\begin{theorem}\label{thm I-III}
    Let $\mathcal{M}$ be a separable factor of type $\mathrm{I}$ or $\mathrm{III}$.
    Then every nonzero operator in $\mathcal{M}$ is the product of two irreducible operators in $\mathcal{M}$.
\end{theorem}

\begin{proof}
    Let $T$ be a nonzero operator in $\mathcal{M}$.
    Suppose that $\mathcal{M}$ is of type $\mathrm{III}$.
    If $T$ is an irreducible operator in $\mathcal{M}$, then we finish the proof by \Cref{lem IR=IR2}.
    If $T$ is reducible, then there exists a projection $P$ in $\mathcal{M}$ such that $P\ne 0,I$, and $PT=TP$.
    Since $I\sim P\sim I-P$, we complete the proof by \Cref{lem half-exact}.

    Suppose that $\mathcal{M}$ is of type $\mathrm{I}_\infty$, i.e., $\mathcal{M}\cong\mathcal{B}(\mathcal{H})$, where $\mathcal{H}$ is a separable infinite-dimensional complex Hilbert space.
    If $W^*(T)'$ is infinite dimensional, then there exists a projection $P$ in $W^*(T)'$ such that both $P$ and $I-P$ are infinite-dimensional projections in $\mathcal{B}(\mathcal{H})$.
    Hence, we complete the proof by \Cref{lem half-exact}.
    If $W^*(T)'$ is finite dimensional, then there are minimal projections $P_1,P_2,\ldots,P_m$ in $W^*(T)'$ such that $I=\sum_{j=1}^{m}P_j$.
    By the minimality of $P_j$, $TP_j$ is  irreducible  in $P_j\mathcal{B}(\mathcal{H})P_j=\mathcal{B}(P_j\mathcal{H})$ for each $1\leqslant j\leqslant m$.
    Without loss of generality, we may assume that $P_1$ is an infinite-dimensional projection and $I-P_1$ is a finite-dimensional projection in $\mathcal{B}(\mathcal{H})$ by \Cref{lem half-exact}.
    It follows that $TP_j$ has closed range for $2\leqslant j\leqslant m$.
    Suppose that $TP_j\ne 0$ for $1\leqslant j\leqslant k$ and $TP_j=0$ for $k+1\leqslant j\leqslant m$.
    Let $P=\sum_{j=1}^{k}P_j$.
    Then $TP$ is the product of a nonzero irreducible operator and an invertible irreducible operator in $\mathcal{B}(P\mathcal{H})$ by \Cref{prop closed-range}.
    Thus, the proof is completed by \Cref{cor add-0s}.
    
    If $\mathcal{M}$ is of type $\mathrm{I}_n$, i.e., $\mathcal{M}\cong M_n(\mathbb{C})$, then $W^*(T)'$ is clearly finite dimensional.
    The proof is a slight modification of the type $\mathrm{I}_\infty$ case.
    This completes the proof.
\end{proof}

\section{Products of irreducible operators in factors of type \texorpdfstring{$\mathrm{II}$}{II}} \label{section-4}

In this section, we will show that every operator is the product of two irreducible operators in separable factors of type $\mathrm{II}$ (see \Cref{thm II-infty} and \Cref{thm II-1}).
The following lemma is a routine construction in factors of type $\mathrm{II}_\infty$.

\begin{lemma}\label{lem sum-finite-projection}
    Let $\mathcal{M}$ be a separable factor of type $\mathrm{II}_\infty$.
    If $\{P_j\}_{j\in\Lambda}$ is an orthogonal family of nonzero finite projections in $\mathcal{M}$ with $I=\sum_{j\in\Lambda}P_j$, then there is a partition $\{\Lambda_1,\Lambda_2\}$ of $\Lambda$ such that $I\sim\sum_{j\in\Lambda_1}P_j\sim\sum_{j\in\Lambda_2}P_j$.
\end{lemma}

\begin{proof}
    Let $\tau$ be a normal faithful semifinite tracial weight on $\mathcal{M}$.
    The separability of $\mathcal{M}$ yields that $\Lambda$ is countable.
    Thus, let $\Lambda$ be the set of positive integers.
    Then
    \begin{equation*}
        \sum_{j=1}^{\infty}\tau(P_j)=\tau(I)=\infty\quad\text{and}\quad
        \tau(P_j)<\infty~\text{for each}~j\geqslant 1.
    \end{equation*}
    Clearly, there exists an integer $n_1>1$ such that $\sum_{j=1}^{n_1}\tau(P_j)\geqslant 1$, and there exists an integer $n_k>n_{k-1}$ such that $\sum_{j=n_{k-1}+1}^{n_k}\tau(P_j)\geqslant 1$ for every $k\geqslant 2$.
    Let
    \begin{equation*}
        \Lambda_1=\{1,2,\ldots,n_1\}\cup
        \bigcup_{k=1}^{\infty}\{n_{2k}+1,n_{2k}+2,\ldots,n_{2k+1}\}\quad\text{and}\quad
        \Lambda_2=\Lambda\setminus\Lambda_1.
    \end{equation*}
    This completes the proof.
\end{proof}

\begin{theorem}\label{thm II-infty}
    Let $\mathcal{M}$ be a separable factor of type $\mathrm{II}_\infty$.
    Then every operator in $\mathcal{M}$ is the product of two irreducible operators in $\mathcal{M}$.
\end{theorem}

\begin{proof}
    Let $T$ be an operator in $\mathcal{M}$ and we denote by $\mathcal{N}$ the relative commutant $W^*(T)'\cap\mathcal{M}$.
    Let $\{P_j\}_{j\in\Lambda}\subseteq\mathcal{N}$ be a maximal orthogonal family of nonzero finite projections and $P=\sum_{j\in\Lambda}P_j$.
    If $P=I$, then there exists a projection $Q$ in $\mathcal{N}$ such that $I\sim Q\sim I-Q$ in $\mathcal{M}$ by \Cref{lem sum-finite-projection}.
    In this case, we complete the proof by \Cref{lem half-exact}.
    Thus, we may assume that $Q=I-P\ne 0$.
    By the maximality of $\{P_j\}_{j\in\Lambda}$, every nonzero projection in $Q\mathcal{N}Q$ is an infinite projection in $\mathcal{M}$.
    In particular, $Q$ is an infinite projection in $\mathcal{M}$.
    If $Q\mathcal{N}Q$ is one dimensional, i.e., $Q\mathcal{N}Q=\mathbb{C}Q$, then $TQ$ is irreducible in $Q\mathcal{M}Q$ and the proof is completed by \Cref{lem half-exceed}.
    Otherwise, there are nonzero projections $Q_1$ and $Q_2$ in $\mathcal{N}$ with $Q_1+Q_2=Q$.
    Then both $Q_1$ and $Q_2$ are infinite projections in $\mathcal{M}$.
    As a consequence, both $Q_1$ and $I-Q_1$ are infinite projections in $\mathcal{M}$.
    We finish the proof by \Cref{lem half-exact}.
\end{proof}

Note that the proof of \Cref{thm II-infty} also works for type $\mathrm{I}_{\infty}$ factors. In terms of \Cref{thm I-III} and \Cref{thm II-infty}, every nonzero operator in $\mathcal{M}$ is the product of two irreducible operators whenever $\mathcal{M}$ is a separable factor of type $\mathrm{I}_n$, $\mathrm{I}_\infty$, $\mathrm{II}_\infty$, or $\mathrm{III}$.
At the end of this paper, we focus on type $\mathrm{II}_1$ factors.

We continue our discussion with the following observation.
Let $\mathcal{M}$ be a separable factor of type $\mathrm{II}_1$ and $T$ a normal operator in $\mathcal{M}$.
Then $T$ lies in a maximal abelian von Neumann subalgebra $\mathcal{A}$ of $\mathcal{M}$.
By \cite[Theorem 5.6.2]{SS-book}, there is a projection $P$ in $\mathcal{A}$ such that $P\sim I-P$.
Thus, $T$ is the product of two irreducible operators in $\mathcal{M}$ by \Cref{lem half-exact}.
Note that every operator in $\mathcal{M}$ is the product of a unitary operator and a positive operator by polar decomposition.
Therefore, every operator in $\mathcal{M}$ is the product of four irreducible operators.
To get a result similar to \Cref{thm I-III} and \Cref{thm II-infty}, we need to develop new techniques in factors of type $\mathrm{II}_1$.

The following lemma is prepared for the proof of \Cref{prop key}.

\begin{lemma}\label{lem partial-isometry}
    Let $\mathcal{M}$ be a separable factor, $\{P_j\}_{j=1}^{\infty}$ a sequence of nonzero orthogonal projections in $\mathcal{M}$, and $P=\sum_{j=1}^{\infty}P_j$.
    Then for any nonzero projection $Q$ in $\mathcal{M}$, there exists a partial isometry $V$ in $P\mathcal{M}Q$ such that $P_jV\ne 0$ for each $j\geqslant 1$.
\end{lemma}

\begin{proof}
    Since $\mathcal{M}$ is a factor, there exists a partial isometry $V_j$ in $P_j\mathcal{M}Q$ such that either $V_j^*V_j=Q$ or $V_jV_j^*=P_j$.
    In particular, $V_j\ne 0$ for each $j\geqslant 1$.
    Write $T := \sum_{j=1}^{\infty}2^{-j}V_j$ in $P\mathcal{M}Q$.
    By polar decomposition, there exists a partial isometry $V$ in $P\mathcal{M}Q$ and a positive operator $H$ in $Q\mathcal{M}Q$ such that $T=VH$.
    Since $P_jVH=P_jT=2^{-j}V_j\ne 0$, we have $P_jV\ne 0$ for each $j\geqslant 1$.
    This completes the proof.
\end{proof}

Let $\mathcal{M}$ be a separable factor, $T$ an operator in $\mathcal{M}$, and $\{P_j\}_{j=1}^k$ finitely many nonzero projections in $\mathcal{M}$ such that $I=\sum_{j=1}^kP_j$ and $P_jT=TP_j$ for each $1\leqslant j\leqslant k$.
Suppose that $TP_j$ is nonzero and irreducible in $P_j\mathcal{M}P_j$ for each $1\leqslant j\leqslant k$.
If $T$ has closed range, then $T$ is the product of two irreducible operators in $\mathcal{M}$ by \Cref{lem IR=IR2} and \Cref{prop closed-range}.
Note that we can remove the conditions $k$ being finite and $T$ having closed range by \Cref{prop key}, which is the main tool in the proof of \Cref{thm II-1}.

\begin{proposition}\label{prop key}
    Let $\mathcal{M}$ be a separable factor and $\{P_j\}_{j=1}^{N}$ nonzero projections in $\mathcal{M}$ such that $I=\sum_{j=1}^{N}P_j$, where $2\leqslant N\leqslant\infty$.
    Suppose that $A$ and $B$ are operators in $\mathcal{M}$ such that $AP_j=P_jA$ and $BP_j=P_jB$ are nonzero irreducible operators in $P_j\mathcal{M}P_j$ for every $j\geqslant 1$.
    If $BP_1$ is invertible in $P_1\mathcal{M}P_1$, $AP_2$ is invertible in $P_2\mathcal{M}P_2$, and $P_1\precsim P_2$, then $AB$ is the product of two irreducible operators in $\mathcal{M}$.
\end{proposition}

\begin{proof}
    We only prove the case $N=\infty$ because the case $2\leqslant N<\infty$ can be proved similarly.
    Let $A_j=AP_j$ and $B_j=BP_j$ for every $j\geqslant 1$.
    Then there exist scalars $\{\lambda_j\}_{j=3}^{\infty}$ in the closed interval $[1,2]$ such that
    \begin{equation*}
        \|\lambda_jA_j\|\ne\|\lambda_kA_k\|\quad\text{and}\quad
        \|\lambda_j^{-1}B_j\|\ne\|\lambda_k^{-1}B_k\|\quad\text{for all}~j\ne k\geqslant 3.
    \end{equation*}
    Since $A_2$ and $B_1$ are invertible, there are scalars $\lambda_2$ large and $\lambda_1$ small such that
    \begin{enumerate} [label= $(\arabic*)$, ref= \arabic*, labelsep = 0.5 em]
        \item  $\sigma_{P_2\mathcal{M}P_2}(\lambda_2A_2)
          \cap\sigma_{P_j\mathcal{M}P_j}(\lambda_jA_j)= \varnothing$ for all $j\ne 2$;
        \item  $\sigma_{P_1\mathcal{M}P_1}(\lambda_1^{-1}B_1)
          \cap\sigma_{P_j\mathcal{M}P_j}(\lambda_j^{-1}B_j)= \varnothing$ for all $j\ne 1$;
        \item  $\|\lambda_jA_j\|\ne\|\lambda_kA_k\|$ and $\|\lambda_j^{-1}B_j\|\ne\|\lambda_k^{-1}B_k\|$ for all $j\ne k$.
    \end{enumerate}
    Let $A'=\sum_{j=1}^{\infty}\lambda_jA_j$ and $B'=\sum_{j=1}^{\infty}\lambda_j^{-1}B_j$.
    Then $AB=A'B'$.
    Without loss of generality, we assume that
    \begin{enumerate}  [label= $(\arabic*')$, ref= $\arabic*'$, labelsep = 0.5 em]
        \item \label{item:b1} $\sigma_{P_2\mathcal{M}P_2}(A_2)
          \cap\sigma_{P_j\mathcal{M}P_j}(A_j)= \varnothing$ for all $j\ne 2$;
        \item \label{item:b2} $\sigma_{P_1\mathcal{M}P_1}(B_1)
          \cap\sigma_{P_j\mathcal{M}P_j}(B_j)= \varnothing$ for all $j\ne 1$;
        \item \label{item:b3} $\|A_j\|\ne\|A_k\|$ and $\|B_j\|\ne\|B_k\|$ for all $j\ne k$.
    \end{enumerate}
    Since $P_1\precsim P_2$, there exists a partial isometry $U$ in $P_2\mathcal{M}P_1$ such that $U^*U=P_1$ and $UU^*\leqslant P_2$.
    Let $P'_3=I-P_1-P_2$, $X_3=A-A_1-A_2$, and $Y_3=B-B_1-B_2$.
    By \Cref{lem partial-isometry}, there exists a partial $V$ in $P'_3\mathcal{M}P_1$ such that $A_jV\ne 0$ for all $j\geqslant 3$.
    We construct two operators in $\mathcal{M}$ as follows
    \begin{equation*}
        X:= A_1+A_2+X_3+A_2U+X_3V\quad\text{and}\quad
        Y:= B_1+B_2+Y_3-UB_1-VB_1.
    \end{equation*}
    For the reader's convenience, we write $X$ and $Y$ in matrix form as follows
    \begin{equation*}
        X=
        \begin{pmatrix}
            A_1 & 0 & 0 \\
            A_2U & A_2 & 0 \\
            X_3V & 0 & X_3
        \end{pmatrix}
        \begin{matrix}
            \operatorname{ran} P_1 \\
            \operatorname{ran} P_2 \\
            \operatorname{ran} P'_3
        \end{matrix}\quad\text{and}\quad
        Y=
        \begin{pmatrix}
            B_1 & 0 & 0 \\
            -UB_1 & B_2 & 0 \\
            -VB_1 & 0 & Y_3
        \end{pmatrix}
        \begin{matrix}
            \operatorname{ran} P_1 \\
            \operatorname{ran} P_2 \\
            \operatorname{ran} P'_3
        \end{matrix}.
    \end{equation*}
    It is clear that $AB=XY$.
    We claim that both $X$ and $Y$ are irreducible in $\mathcal{M}$.
    
    Let $Q$ be a projection in $\mathcal{M}$.
    For convenience, let $Q_{ij}=P_iQP_j$ for $i,j\geqslant 1$.
    If $Q$ commutes with $X$, then $Q_{2j}=0$ for all $j\ne 2$ by the condition \eqref{item:b1} and \Cref{lem rosenblum}.
    Hence $Q_{22}$ is a projection in $P_2\mathcal{M}P_2$ commuting with $A_2$.
    By the irreducibility of $A_2$, we may assume that $Q_{22}=0$.
    Then $A_2UQ_{11}=0$.
    Since $A_2$ is invertible and $U^*U=P_1$, we have $Q_{11}=0$.
    It follows that $Q_{1j}=0$ for all $j\geqslant 1$.
    Thus, $Q\leqslant P'_3$ and $QX_3=X_3Q$.
    By the condition \eqref{item:b3} and \Cref{lem IR-sum-commutant}, we have $Q_{ij}=0$ for $i\ne j$, and $Q_{jj}=0$ or $P_j$ for $j\geqslant 3$.
    Since $QX_3V=0$, we have $Q_{jj}A_jV=0$.
    Hence $Q_{jj}=0$ for each $j\geqslant 3$.
    Therefore, $Q=0$ and $X$ is irreducible in $\mathcal{M}$.
    
    If $Q$ commutes with $Y$, then $Q_{1j}=0$ for all $j\ne 1$ by the condition \eqref{item:b2} and \Cref{lem rosenblum}.
    Hence $Q_{11}$ is a projection in $P_1\mathcal{M}P_1$ commuting with $B_1$.
    By the irreducibility of $B_1$, we may assume that $Q_{11}=0$.
    It follows that $Q\leqslant P_2+P'_3$ and $Q(B_2+Y_3)=(B_2+Y_3)Q$.
    By the condition \eqref{item:b3} and \Cref{lem IR-sum-commutant}, we have $Q_{ij}=0$ for $i\ne j$, and $Q_{jj}=0$ or $P_j$ for $j\geqslant 2$.
    It follows that $Q_{22}UB_1=0$ and $QP'_3VB_1=0$.
    Then $Q_{22}=0$ and $Q_{jj}P_jV=0$ for each $j\geqslant 3$.
    Hence $Q_{jj}=0$ for each $j\geqslant 3$.
    Therefore, $Q=0$ and $Y$ is irreducible in $\mathcal{M}$.
    We complete the proof.
\end{proof}

Recall that the {\em central support $C_P$} of a projection $P$ in a von Neumann algebra is the minimal central projection with the property $C_PP=P$.

\begin{lemma}\label{lem atomic}
    Let $\mathcal{N}$ be a von Neumann algebra, $\{P_j\}_{j\in\Lambda}$ a maximal orthogonal family of minimal projections in $\mathcal{N}$, and $P=\sum_{j\in\Lambda}P_j$.
    Then $P$ is a central projection in $\mathcal{N}$.
    Moreover, $\mathcal{N}(I-P)$ is diffuse or $I-P=0$.
\end{lemma}

\begin{proof}
    If $P$ is not a central projection in $\mathcal{N}$, then $C_PC_{I-P}\geqslant C_P(I-P)=C_P-P\ne 0$.
    It follows that $C_{P_j}C_{I-P}\ne 0$ for some $j\in\Lambda$.
    By \cite[Proposition 6.1.8]{KR2}, there is a projection $P_0$ in $\mathcal{N}$ such that $P_j\sim P_0\leqslant I-P$.
    Then $P_0$ is also a minimal projection in $\mathcal{N}$.
    That contradicts the maximality of $\{P_j\}_{j\in\Lambda}$.
    Therefore, $P$ is a central projection in $\mathcal{N}$.
    In particular, if $I-P\ne 0$, then $\mathcal{N}(I-P)$ is diffuse by the maximality of $\{P_j\}_{j\in\Lambda}$.
    This completes the proof.
\end{proof}

Note that for each finite factor $\mathcal{M}$, there is a unique normal faithful tracial state $\tau$ such that $P\precsim Q$ if and only if $\tau(P)\leqslant\tau(Q)$ for all projections $P$ and $Q$ in $\mathcal{M}$.

\begin{theorem}\label{thm II-1}
    Let $\mathcal{M}$ be a separable factor of type $\mathrm{II}_1$.
    Then every operator in $\mathcal{M}$ is the product of two irreducible operators in $\mathcal{M}$.
\end{theorem}

\begin{proof}
    Let $\tau$ be the unique normal faithful tracial state on $\mathcal{M}$.
    Let $T$ be an operator in $\mathcal{M}$, $\mathcal{N}$ the relative commutant $W^*(T)'\cap\mathcal{M}$, $\{P_j\}_{j\in\Lambda}$ a maximal orthogonal family of minimal projections in $\mathcal{N}$, and $P=\sum_{j\in\Lambda}P_j$.
    Then $P$ is a central projection by \Cref{lem atomic}.
    Moreover, if $I-P$ is nonzero, then $\mathcal{N}(I-P)$ is diffuse.
    Clearly, there are projections $Q_1$ and $Q_2$ in $\mathcal{N}(I-P)$ such that $I-P=Q_1+Q_2$ and $\tau(Q_1)=\tau(Q_2)$.
    Thus, $T(I-P)$ is the product of two irreducible operators in $(I-P)\mathcal{M}(I-P)$ by \Cref{lem half-exact}.
    The proof is divided into three cases.
    
    \noindent \textbf{Case I.}
    $\Lambda= \varnothing$.
    Then $P=0$, $\mathcal{N}$ is diffuse, and we finish the proof.
    
    \noindent \textbf{Case II.}
    $|\Lambda|=1$.
    Then $TP$ is irreducible in $P\mathcal{M}P$.
    If $\tau(P)\geqslant\frac{1}{2}$, then we complete the proof by \Cref{lem half-exceed}.
    If $\tau(P)<\frac{1}{2}$, then there exists a projection $Q$ in $\mathcal{N}(I-P)$ with $\tau(Q)=\frac{1}{2}$.
    We complete the proof by \Cref{lem half-exact}.
    
    \noindent \textbf{Case III.}
    $|\Lambda|\geqslant 2$.
    Without loss of generality, we assume that $\{P_j\}_{j\in\Lambda}=\{P_j\}_{j=1}^N$ and $\tau(P_1)\leqslant\tau(P_2)$, where $2\leqslant N\leqslant\infty$.
    Then $TP_j$ is irreducible in $P_j\mathcal{M}P_j$ for each $j$.
    Moreover, $T(I-P)$ is the product of two irreducible operators in $(I-P)\mathcal{M}(I-P)$ or $I-P=0$.
    By \Cref{lem IR=IR2} and \Cref{prop key}, $T$ is the product of two irreducible operators in $\mathcal{M}$.
\end{proof}

\begin{remark}\label{rem II-1-non-Gamma}
    The separability of $\mathcal{M}$ ensures the existence of irreducible operators in the reduced von Neumann algebra $P\mathcal{M}P$ for each nonzero projection $P$ in $\mathcal{M}$ in the proof of \Cref{lem half-exact}.
    By \cite[Theorem 6.6]{SS19} and the proof of \Cref{thm II-1}, every operator in nonseparable non-$\Gamma$ type $\mathrm{II}_1$ factors is the product of two irreducible operators.
\end{remark}


\begin{thebibliography}{99}

\bibitem{CFY21}
Xinyan Cao, Junsheng Fang, Zhaolin Yao.
{\em Strong sums of projections in type  $\mathrm{II}_1$  factors.}
J. Funct. Anal. 281 (2021), no. 5, Paper No. 109088, 11 pp.
MR4253932

\bibitem{CFY24}
Xinyan Cao, Junsheng Fang, Zhaolin Yao.
{\em On finite sums of projections and Dixmier's averaging theorem for type  $\mathrm{II}_1$  factors.}
J. Funct. Anal. 287 (2024), no. 8, Paper No. 110568, 29 pp.
MR4777788


\bibitem{FT}
P. Fillmore, D. Topping.
{\em Sums of irreducible operators.}
Proc. Amer. Math. Soc. 20 (1969), 131-133.
MR0233226


\bibitem{HMS26}
Donald Hadwin, Minghui Ma, Junhao Shen, 
{\em Voiculescu's theorem in properly infinite factors.}
J. Funct. Anal. 290 (2026), no. 1, Paper No. 111198, 34 pp.
MR4964138


\bibitem{Hal}
P. Halmos.
{\em Irreducible operators.}
Michigan Math. J. 15 (1968), 215-223.
MR0231233


\bibitem{Her89-book}
D. Herrero.
{\em  Approximation of Hilbert space operators. Vol. 1.}
Pitman Res. Notes Math. Ser., 224,
Longman Scientific \& Technical, Harlow; copublished in the United States with John Wiley \& Sons, Inc., New York, 1989, xii+332 pp.
MR1088255


\bibitem{JW98-book}
C. Jiang, Z. Wang.
{\em Strongly irreducible operators on Hilbert space.}
Pitman Res. Notes Math. Ser., 389, 
Longman, Harlow, 1998. x+243 pp.
MR1640067


\bibitem{JW06}
C. Jiang, Z. Wang.
{\em Structure of Hilbert space operators.}
World Scientific Publishing Co. Pte. Ltd., Hackensack, NJ, 2006. x+248 pp.
MR2221863


\bibitem{JW98}
C. Jiang, P. Wu.
{\em Sums of strongly irreducible operators.}
Houston J. Math. 24 (1998), no. 3, 467-481.
MR1686618



\bibitem{KR2}
R.\,Kadison, J.\,Ringrose.
{\em Fundamentals of the theory of operator algebras, II, Advanced theory.}
Academic Press, Orlando, FL, 1986.
MR0859186


\bibitem{Rad}
H. Radjavi.
{\em Every operator is the sum of two irreducible ones.}
Proc. Amer. Math. Soc. 21 (1969), 251-252.
MR0238112


\bibitem{RR}
H. Radjavi, P. Rosenthal.
{\em The set of irreducible operators is dense.}
Proc. Amer. Math. Soc. 21 (1969), 256.
MR0236748


\bibitem{RR73}
H. Radjavi, P. Rosenthal.
{\em Invariant subspaces.}
Ergeb. Math. Grenzgeb., Band 77
[Results in Mathematics and Related Areas]
Springer-Verlag, New York-Heidelberg, 1973. xi+219 pp.
MR0367682


\bibitem{Ros}
M. Rosenblum.
{\em On the operator equation $BX-XA=Q$.}
Duke Math. J. 23 (1956), 263-269.
MR0079235


\bibitem{SS19}
J. Shen, R. Shi.
{\em Reducible operators in non-$\Gamma$ type $\mathrm{II}_1$ factors.}
arXiv:1907.00573v3, 2019.


\bibitem{SS20}
J. Shen, R. Shi.
{\em Sum of irreducible operators in von Neumann factors.}
Proc. Amer. Math. Soc. 148 (2020), no. 7, 2901-2908.
MR4099778

\bibitem{SS-book}
A. Sinclair, R. Smith.
{\em Finite von Neumann algebras and masas.}
London Math. Soc. Lecture Note Ser., 351 Cambridge University Press, Cambridge, 2008. x+400 pp.
MR2433341

\bibitem{WFS}
S. Wen, J. Fang, R. Shi.
{\em On irreducible operators in factor von Neumann algebras.}
Linear Algebra Appl. 565 (2019), 239-243.
MR3892720


\end{thebibliography}
\end{document}